\documentclass{amsart}
\newtheorem{Theo}{Theorem}
\newtheorem{Lem}{Lemma}
\newtheorem{Prop}{Proposition}
\newtheorem{Cor}{Corollary}
\newtheorem{Problem}{Problem}
\newcommand{\rg}{\mathrm{rg}}
\newcommand{\Hom}{\mathrm{Hom}}
\newcommand{\de}{\mathrm{def}_p}
\begin{document}
\title{A $p$-group with positive Rank Gradient}
\author[J.-C. Schlage-Puchta]{Jan-Christoph Schlage-Puchta}
\begin{abstract}
We construct for $d\geq 2$ and $\epsilon>0$
a $d$-generated $p$-group $\Gamma$, which in an asymptotic sense
behaves almost like a $d$-generated free pro-$p$-group. We show
that a subgroup of index $p^n$ needs $(d-\epsilon)p^n$ generators, and
that the subgroup growth of $\Gamma$ satisfies
$s_{p^n}(\Gamma)>s_{p^n}(F_d^p)^{1-\epsilon}$, where $F_d^p$ is the
  $d$-generated free pro-$p$-group. To do so we introduce a new
invariant for finitely generated groups and study some of its basic
properties. 
\end{abstract}
\maketitle
MSC-Index: 20F05, 20F50, 20F69 
\section{Introduction and results}
Burnside asked in 1902 whether a finitely generated torsion group is
necessary finite. This question remained open until 1964, when
Golod\cite{Golod} 
constructed an infinite, finitely generated
$p$-group. In the mean time there have been several other
constructions of finitely generated infinite $p$-groups, we only note
the construction of Gupta and Sidki\cite{GuSi} using 
automorphisms of the $p$-adic tree. After
Burnside's question was settled, most further investigations considered
the question whether finitely generated periodic could be replaced by
a stronger restriction which still allows for infinite groups. Much
less work was done concerning the question as to how large a finitely
generated torsion group could become. Here we show that such a group can in
some sense be almost as large as a free group with the same number of
generators. Denote by $d(G)$ the number of generators of $G$. Let
$\Gamma$ be a finitely generated residually finite infinite group,
$\mathcal{N}=(N_i)_{i=1}^\infty$ be a sequence of finite index normal
subgroups of $\Gamma$ with $\bigcap N_i = \{1\}$. Define the rank
gradient $\rg(\Gamma, \mathcal{N})$ of $\Gamma$ along $\mathcal{N}$ as
$\liminf\frac{d(N_i)}{(\Gamma:N_i)}$, and the rank gradient of
$\Gamma$ as the infimum of $\rg(\Gamma, \mathcal{N})$ taken over all
sequences $N_i$ of finite index normal subgroups with trivial
intersection. Obviously, $\rg(\Gamma)\leq d(\Gamma)-1$. This invariant
was introduced by Lackenby\cite{Lack} in connection 
with the virtual Haken conjecture for 3-manifolds, Ab\'ert and
Nikolov\cite{AbNi} showed that it is also related to topological
dynamics. In relation with the virtual positive $b_1$-conjecture,
Nikolov asked whether there exists a group of positive rank gradient,
such that every finite index subgroup does not map onto an infinite
cyclic group. Here we show that such a group indeed exists, in fact,
we prove the following.
\begin{Theo}
\label{thm:Rankgrad}
Let $d\geq 2$ be an integer, $p$ a prime, and $\epsilon>0$. Then there
exists a $d$-generated $p$-group $\Gamma$ with $\rg(\Gamma)\geq
d-1-\epsilon$.
\end{Theo}

The existence of finitely generated $p$-groups with positive rank
gradient was independently shown by Osin\cite{Osin}. Using results by
Lackenby\cite{LackAm} and Abert, Jaikin-Zapirai and Nikolov\cite{AJN},
we obtain the following. 
\begin{Cor}
\label{Cor}
There exists an infinite non-amenable torsion group.
\end{Cor}

This result was recently proven by Ershov\cite{Er} (see also
\cite[Corollary~1.3]{Osin}), however, the proof given here is more
elementary. 

For the proof we define another invariant, the $p$-deficiency of a
group $\Gamma$. Let $\Gamma$ be a group with a presentation $\langle
X|R\rangle$, where we assume $X$ to be finite. We can view an element
$r\in R$ as an element in the free group $F(X)$ over $X$, and define
the $p$-order $\nu_p(r)$ of $r$ to be the largest integer $k$,
such that there is 
some $s\in F(X)$ with $r=s^{p^k}$. Define the $p$-deficiency of
$\langle X|R\rangle$ as $|X|-1-\sum_{r\in R}p^{-\nu_p(r)}$, and the
$p$-deficiency $\de(\Gamma)$ of $\Gamma$ as the supremum taken over all
representations of $\Gamma$ with finite generating set. Our technical
main result is the following.

\begin{Theo}
\label{thm:Multi}
Let $\Gamma$ be an infinite group, $\Delta\triangleleft\Gamma$ be a
normal subgroup of index $p$. Then $\de(\Delta)\geq
p\cdot\de(\Gamma)$.
\end{Theo}

From this we deduce the following. For a group $\Gamma$ denote by
$s_n(\Gamma)$ the number of subgroups of index $n$, and by
$s_n^{\triangleleft\triangleleft}(\Gamma)$ the number of subnormal
  subgroups of index $n$.

\begin{Prop}
\label{prop:list}
\begin{enumerate}
\item If $\Gamma$ is a $p$-group which is residually finite, then
  $\de(\Gamma)\leq\rg(\Gamma)$.
\item For every $\epsilon>0$ we have $s_{p^n}^{\triangleleft\triangleleft}
  (\Gamma)\gg p^{(\de(\Gamma)-\epsilon)\frac{p^{n+1}}{p-1}}$.
\item For every $d\geq 2$, $p$ a prime and $\epsilon>0$ there exists a
  $d$-generated $p$-group $\Gamma$ with $\de(\Gamma)\geq
  d-1-\epsilon$.
\item For every $\epsilon>0$ there exists a $d$-generated
  pro-$p$-group $\Gamma$ with $s_{p^n}(\Gamma)\gg
  s_{p^n}(\widehat{F_d})^{1-\epsilon}$, where $\widehat{F_d}$ is the
  $d$-generated free pro-$p$-group.
\end{enumerate}
\end{Prop}

We can define the $p$-deficiency for pro-finite groups as we did for
discrete groups, we only have to replace the rank by the topological
rank, that is, the minimal number of elements generating a dense
subgroup. In the case of pro-$p$-groups this 
allows us to construct an Euler-characteristic for finitely generated
pro-$p$-groups as follows. We define
\[
\chi_p(\Gamma) = -\sup_{\Delta<\Gamma, \Delta\mbox{\footnotesize\ open}}
\frac{\de(\Delta)}{(\Gamma:\Delta)}.
\]
Note that the sign is somewhat arbitrary, we chose it to be consistent
with the Euler-characteristic of a virtually free group.
For this invariant we prove the following.
\begin{Prop}
\label{Prop:Euler}
\begin{enumerate}
\item If $\Gamma$ is a pro-$p$-group, then $-\rg(\Gamma)\leq\chi_p(\Gamma)\leq-\de(\Gamma)$.
\item If $\chi_p(\Gamma)<0$, and $\Delta$ is an open subgroup
  of $\Gamma$, then $\chi_p(\Delta)=(\Gamma:\Delta)\chi_p(\Gamma)$.
\item Let $\Gamma$ be a virtually free group, which contains a normal free
  subgroup of $p$-power index, $\widehat{\Gamma}$ be the
  pro-$p$-completion of $\Gamma$. Then $\chi_p(\widehat{\Gamma})=\chi(\Gamma)$.
\item Let $\Gamma$ be a Fuchsian group, which contains a normal surface group
  as a subgroup of $p$-power index, $\widehat{\Gamma}$ be the
  pro-$p$-completion of $\Gamma$. Then $\chi_p(\widehat{\Gamma})=-\mu(\Gamma)$,
  where $\mu$ is the hyperbolic volume of a group.
\end{enumerate}
\end{Prop}
In spite of Theorem~\ref{thm:Rankgrad}, positive $p$-deficiency
implies that a group behaves roughly like a large group. For example,
the following is an immediate consequence of Lackenby's
characterisation of finitely presented large groups (confer \cite{LackLarge}).

\begin{Prop}
\label{prop:large}
A finitely presented group $\Gamma$ with positive $p$-deficiency contains a
finite index subgroup $\Delta$ which projects surjectively onto a
non-abelian free group.
\end{Prop}

I would like to thank the referee for pointing me to the work of Osin
and Ershov, and M. Ab\'ert for simplifying the proof of
Theorem~\ref{thm:Multi}.

\section{Proof of Theorem 2}

Let $\Gamma$ be a group with a presentation $\langle X|R\rangle$,
where $|X|=d$ is finite. Let $F$ be the free group over $X$. For a subgroup
$H\leq F$ and an element $g\in F$, we define $g^H=\{g^h:h\in H\}$ to
be the $H$-conjugacy class of $g$. Subgroups of free groups are free,
since the $p$-order of an element is only defined with respect to one
fixed free group containing it we define $\nu_{p, H}(g)$ for $g\in H$ to be the
largest $k$, such that there exists some $h\in H$ with $h^{p^k}=g$. 
\begin{Lem}
\label{Lem:trivial}
Suppose that $N\triangleleft F$ is a normal subgroup of index
$p$. Then $\nu_{p, N}(g)\geq\nu_{p, F}(g)-1$ for all $g\in N$.
\end{Lem}
\begin{proof}
If $h^{p^{k}}=g$ then $h^p\in N$ and so $\nu_{p, N}(g)\geq k-1$.
\end{proof}
\begin{Lem}
\label{Lem:split}
Let $N\vartriangleleft F$ be of index $p$ and let $g\in N$. Then we
have either $g^{F}=g^{N}$ or  
\begin{equation*}
g^{F}=\bigcup\limits_{i=1}^{p}g_{i}^{N}.
\end{equation*}
Moreover, in the latter case we have $\nu_{p, N}(g)\leq\nu_{p, F}(g)$.
\end{Lem}
\begin{proof}
The group $F$ acts on $g^{F}$ transitively by
conjugation, so the normal subgroup $N$ has either $1$ or $p$ orbits. The
latter is equivalent to saying that the stabilizer of $g$ in $H$ equals the
stabilizer of $g$ in $F$, that is, that the centralizer $C_{F}(g)\leq N$.
Assume that this is the case. If $h^{p^{k}}=g$ then $h$ and $g$ commute,
hence $h\in N$. So $\left\vert g\right\vert _{N}=\left\vert g\right\vert
_{F} $.
\end{proof}
Let $\Gamma =\langle X\mid R\rangle$ be a presentation for
$\Gamma$ and let $F$ be the free group over $X$. Let $\phi :F\rightarrow 
\Gamma$ be the homomorphism defined by the presentation, $K$ be the
kernel of $\phi$, and $N=\phi ^{-1}(\Delta)$ be the preimage of $\Delta$. Then
$N\vartriangleleft F$ is a normal subgroup of index $p$ and by the
Nielsen-Schreier theorem $H$ is generated by $n=(d-1)p+1$ elements, we
fix a generating set $Y=y_1,\ldots ,y_n$.

We construct a presentation $\Delta=\langle Y\mid S\rangle$ for $\Delta$ using
this generating set and $K$. Now $K$ is generated by the conjugacy
classes $r^F$, $r\in R$ as a subgroup. 
Let us use Lemma \ref{Lem:split}. If for some $r\in R$ we have
$r^F=r^N$ then we add $r$ to $S$, expressed as a word over $Y$.
In this case, by Lemma \ref{Lem:trivial} we have $\nu_{p, N}(r)\geq
\nu_{p, F}(r)-1$. Otherwise, $r^F$ is the
disjoint union of $p$ conjugacy classes under $N$; let us add one element
from each conjugacy class to the presentation as a relation $s$. In this
case, by Lemma \ref{Lem:split} we have $\nu_{p,N}(r)=\nu_{p,F}(r)$. We
do so for every $r\in R$, and collect the resulting relations into the
set $S$. Then $\langle Y|S\rangle$ will be a presentation for
$\Delta$. Using these estimates we get  
\begin{multline*}
\de\left( \langle Y|S\rangle \right)\geq
n-1-\sum\limits_{s\in S}p^{-\nu_{p, H}(s)}\geq
p(d-1)-p\sum\limits_{r\in R}p^{-\nu_{p, H}(r)} =
\mathrm{def}_{p}\left( \langle X|R\rangle \right).
\end{multline*}
Let $\epsilon>0$ be given. Then there exists a presentation
$\Gamma=\langle X|R\rangle$ such that $\de(\Gamma)\leq\de(\langle
X|R\rangle)-\epsilon$. From the computation above we see that there
exists a presentation $\langle Y|S\rangle$ of $\Delta$ with $\de(\langle
Y|S\rangle)\geq p\cdot\de(\langle X|R\rangle)$, hence, $\de(\Delta)\geq
p\cdot(\de(\Gamma)-\epsilon)$. Since $\epsilon$ is arbitrary,
Theorem~\ref{thm:Multi} follows.

Note that elements of $S$ might be redundant, or there could
be more economic presentations of $\Delta$, hence, in general we do
not have equality.

\section{Proofs of the other statements}

We begin by proving Proposition~\ref{prop:list} (1).
Let $\Delta$ be a normal subgroup of index $p^k$. Then there exists a chain
\[
\Gamma=\Delta_0>\Delta_1>\dots>\Delta_k=\Delta,
\]
such that $(\Delta_i:\Delta_{i+1})=p$. From Theorem~\ref{thm:Multi} we
obtain $\de(\Delta)\geq p^k\de(\Gamma)$. Hence, $\Delta$ has a
presentation $\langle X|R\rangle$, such that $|X|-1-\sum_{r\in R}
p^{-\nu_2(r)}\geq p^k\de(\Gamma)-\epsilon$. Let $R_1\subseteq R$ be the set of
relations, which are not $p$-th powers. Every $r\in R_1$ contributes 1 to
the left hand side sum, hence, 
$|X|-1-|R_1|\geq p^k\de(\Gamma)-\epsilon$.
We now bound $|\Hom(\Delta, C_p)|$. This number equals the number of
solutions of the system $\{r=1|r\in R\}$, where the variables are from
$C_p$, and the equations are to be interpreted as equations in
$C_p$. Since every $p$-th power is trivial in $C_p$, the equations from
$R\setminus R_1$ are trivially satisfied, and we see that the original
system is equivalent to the system
$\{r=1|r\in R_1\}$. This system can be viewed as a system of $|R_1|$
linear equations in $|X|$ variables over the field with $p$ elements,
hence there are at least $p^{|X|-|R_1|}>p^{p^k\de(\Gamma)}$
solutions. Then the $p$-Frattini quotient of $\Delta$ has cardinality
at least $p^{p^k\de(\Gamma)}$, which implies that $\Delta$ has an
elementary abelian p-group of rank at least $p^k\de(\Gamma)$ as
quotient. But then $\Delta$ itself has rank at least
$p^k\de(\Gamma)$. Hence, for every normal subgroup $\Delta$ of
$p$-power index we have that $\Delta$ has rank at least
$(\Gamma:\Delta)\de(\Gamma)$, that is, $\rg(\Gamma)\geq\de(\Gamma)$.

Now we prove part (2). We have
$s_{p^n}^{\triangleleft\triangleleft}(\Gamma) =
s_{p^n}(\widehat{\Gamma})$, where $\widehat{\Gamma}$ is the
pro-$p$-completion of $\Gamma$. For a pro-$p$-group $\Gamma$ let
$\mathrm{sc}_n(\Gamma)$ to be the number of subgroup chains
$\Gamma=\Delta_0>\Delta_1>\dots>\Delta_n$, where
$(\Delta_i:\Delta_{i+1})=p$. The quantities $\mathrm{sc}_n(\Gamma)$ and
$s_{p^n}(\Gamma)$ are linked via the following (confer \cite{large}).
\begin{Lem}
Foir a pro-$p$-group $\Gamma$ we have
\[
\mathrm{sc}_n(\Gamma)\geq s_{p^n}(\Gamma)\gg p^{-n^2}\mathrm{sc}_n(\Gamma).
\]
\end{Lem}
As in the proof of part (1)
we find that the $p$-Frattiniquotient of a subgroup of index $p^k$ in
$\Gamma^p$ has rank at least $p^k\de(\Gamma)$, that is, a subgroup of
index $p^k$ has at least $\frac{p^{p^k}\de(\Gamma)-1}{p-1}$ subgroups of
index $p$. Hence, the number of subgroup chains of length $n$ is at
least
\begin{eqnarray*}
\mathrm{sc}_n(\Gamma^p) & \geq & \frac{p^{p^{k_0}}\de(\Gamma)-1}{p-1}
\frac{p^{p^{k_0+1}}\de(\Gamma)-1}{p-1} \cdots
\frac{p^{p^{n-1}}\de(\Gamma)-1}{p-1}\\
 & \gg & p^{\frac{p^n-1}{(p-1)}\de(\Gamma)} (p-1)^{-n}\\
 & \geq & p^{\frac{p^n-1}{(p-1)}(\de(\Gamma)-\epsilon)},
\end{eqnarray*}
where $k_0$ is chosen sufficiently large to ensure that
$p^{p^{k_0}}\de(\Gamma)-1$ is positive. Our claim now follows by
combining these estimates.

We construct a $d$-generated $p$-group with $p$-deficiency close to
$d-1$ as follows. Number the words in $F_d$ in some way as $w_1,
\ldots$, and choose some integer $k$. Then define the group $\Gamma =
\langle x_1, \ldots, x_d|w_1^{p^k}, w_2^{p^{k+1}}, w_3^{p^{k+2}},
\ldots\rangle$. Every element of $\Gamma$ is a word in $x_1, \ldots,
x_d$, hence, every element has order a power of $p$. The
$p$-deficiency of this presentation is
\[
d-1-\sum_{\nu=k}^\infty = d-1-\frac{p^{1-k}}{p-1} \geq d-1-2p^{-k},
\]
hence, for $k$ large enough we have $\de(\Gamma)>d-1-\epsilon$.

Finally the fourth statement follows immediatelly from the previous two.

Theorem~\ref{thm:Rankgrad} does not follow immediately from
Proposition~\ref{prop:list} (1) and (3), since the group constructed
in (3) is not necessarily residually finite. In fact, factoring
out the residual might change the $p$-deficiency of a group in an
unforeseen way. Let $\Gamma$ be a $p$-group with residual $N$,
and $\Delta$ a finite index subgroup. Then $\de(\Delta)$ is a lower
bound for $\dim\Delta/([\Delta, \Delta]\Delta^p)$, and since $[\Delta,
\Delta]\Delta^p$ contains $N$, $\dim\Delta/([\Delta,
\Delta]\Delta^p)$ is a lower bound for $d(G)$. Hence, when
constructing a group as in (3), its maximal 
residual finite quotient will satisfy the requirements of
Theorem~\ref{thm:Rankgrad}.

Corollary~\ref{Cor} follows immediatelly from
Theorem~\ref{thm:Rankgrad} and the fact that a group of positive rank
gradient cannot be amenable. This was proven by
Lackenby\cite{LackAm} for finitely presented groups, and by Abert,
Jaikin-Zapirai and Nikolov\cite{AJN} in the general case.

Lackenby\cite{LackLarge} showed that a finitely presented group $\Gamma$
which has positive rank gradient along a series $N_1>N_2>\dots$ of
normal subgroups, such that $N_i/N_{i+1}$ is abelian is large. From
this Proposition~\ref{prop:large} follows immediatelly on observing
Proposition~\ref{prop:list}(1).

\section{Towards an Euler characteristic}

Here we prove Proposition~\ref{Prop:Euler}. The inequality
$\chi_p(\Gamma)\leq-\de(\Gamma)$ follows from the fact that the
supremum of a function taken over all finite index subgroups is at
least the value at the group itself. For the inequality
$-\rg(\Gamma)\leq\chi_p(\Gamma)$ note first that since $\Gamma$ is a
pro-$p$-group we have that for all open
subgroups $\Delta$ of $\Gamma$ there exists a chain
$\Gamma=N_0\triangleright N_1\triangleright\dots\triangleright
N_k=\Delta$ with $N_i/N_{i+1}\cong C_p$. Hence we have that
$\de(\Delta)\geq(\Gamma:\Delta)\de(\Gamma)$ for all open subgroups
$\Delta$ of $\Gamma$. 

Let $(N_i)$ be a descending sequence of normal open
subgroups with trivial intersection, such that
$\lim_{i\rightarrow\infty}\frac{d(N_i)}{(\Gamma:N_i)}\leq
\rg(\Gamma)+\epsilon$, and choose an index $i$, such that $d(N_i)\leq
\rg(\Gamma)+2\epsilon$. Let $\Delta$ be an open subgroup of $\Gamma$ with
$\de(\Delta)>(-\chi_p(\Gamma)-\epsilon)(\Gamma:\Delta)$. Then
$\de(\Delta\cap N_i)\geq(\Delta:\Delta\cap N_i)\de(\Delta)$, hence, we
may assume that $\Delta$ is a subgroup of $N_i$. On the
other hand we have $d(\Delta)\leq(N_i:\Delta) d(N_i)$, hence, 
\begin{multline*}
-\chi_p(\Gamma)\leq \frac{\de(\Delta)}{(\Gamma:\Delta)}+\epsilon \leq
\frac{d(\Delta)}{(\Gamma:\Delta)}+\epsilon \leq
\frac{d(N_i)(N_i:\Delta)}{(\Gamma:\Delta)}+\epsilon\\
 \leq \frac{(\rg(\Gamma)+2\epsilon)(\Gamma:N_i)(N_i:\Delta)}
{(\Gamma:\Delta)}+\epsilon = \rg(\Gamma) + 3\epsilon,
\end{multline*}
and our claim follows.

For the multiplicativity note that for $\Gamma>\Delta$ we have
$\chi_p(\Gamma)\geq(\Gamma:\Delta)\chi_p(\Delta)$, since the set over
which the supremum is taken to compute the left hand side is a
superset of the set used for the right hand side. Let $U<\Gamma$ be a
finite index subgroup. Then we
have $\de(U\cap\Delta)\geq(U:U\cap\Delta)\de(U)$, hence
\[
-\chi_p(\Delta)\geq \frac{\de(U\cap\Delta)}{(\Delta:U\cap\Delta)} \geq
\frac{\de(U)(U:U\cap\Delta)}{(\Delta:U\cap\Delta)} =
\frac{\de(U)(\Gamma:\Delta)}{(\Gamma:U)}. 
\]
Now choose for $U$ a subgroup with $\de(U)\geq
(\Gamma:U)(-\chi_p(\Gamma)-\epsilon)$ we obtain $-\chi_p(\Delta)\geq
(\Gamma:\Delta)(-\chi_p(\Gamma)-\epsilon)$, which implies our claim.

If $\Gamma$ is free pro-$p$-group with $d$ generators, then the
obvious presentation 
yields $\de(\Gamma)\geq d-1$. On the other hand the rank gradient of
$\Gamma$ is $d-1$, hence, we have
$\de(\Gamma)=-\chi_p(\Gamma)=d-1$. In the same way we find that if
$\Gamma$ is a $d$-generated Demuskin-group we have
$\de(\Gamma)=-\chi_p(\Gamma)=d-2$. Since the pro-$p$-completion of a
virtually free group, which contains a free normal subgroup of
$p$-power index has normal subgroup of the same index, which is free
pro-$p$, the third statement follows. Similarly, the
pro-$p$-completion of a Fuchsian 
group which contains only elliptic elements of $p$-power order is
virtually free or virtually Demuskin. From this our last claim follows.

We remark that $\chi_p$ is not multiplicative with respect to direct
products. For example, take $\Gamma=\widehat{F_2\times F_2}$.
The rank gradient of this group is the same as the rank gradient of
the discrete group $F_2\times F_2$, and taking normal subgroups of the
form $N\times N$ one sees that the rank gradient of $F_2\times F_2$ is 0.
On the other hand,  $\chi_p(\widehat{F_2})=1$, thus, $\chi_p$ is not
multiplicative on direct products.

\section{Problems}

The most obvious question is as to whether our definition of
$p$-deficiency is the right one. Right now the definition appears
rather ad hoc, and it might be difficult to prove anything about this
invariant without a more conceptual definition. Therefore we pose the
following.

\begin{Problem}
Give a homological characterization of $\de(\Gamma)$.
\end{Problem}

In our definition of $\de$ we weighted a relator $r$ according to the
largest $k$ such that $r=0$ holds true in every group of exponent
$p^k$. However, it appears more natural to weight according to the
largest $k$ such that $r=0$ holds true in $P_k$, the $p$-Sylow
subgroup of the symmetric group of $S_{p^k}$. However, the na\"ive
approach does not work, since all finite index subgroups of the group
$\Gamma=\langle x, y|[x, y]=1\rangle$ are isomorphic to $\Gamma$, that
is, if we want a statement like Theorem~\ref{thm:Multi} to hold, we
have to assign the weight 1 to $[x, y]$, although applying a
commutator pushes elements down the upper central series. Hence, we
are looking for the correct weights.

\begin{Problem}
Find a function $w:F_r\rightarrow[0, 1]$, such that
$w(x^p)\leq\frac{w(x)}{p}$, the $w$-deficiency defined via
$\mathrm{def}_w(\langle X|R\rangle)=|X|-1-\sum_{r\in R} w(r)$ satisfies
Theorem~\ref{thm:Multi}, and $w(x_i)\rightarrow 0$ as $x_i\rightarrow
1$ in the pro-$p$-completion $F_r^p$ of $F_r$.
\end{Problem}

It appears difficult to compute the $p$-deficiency. Giving one example
of a presentation
one can give a lower bound, however, proving the non-existence of a
presentation of a certain form appears difficult, unless one has a
matching upper bound e.g. by subgroup-growth. However, we believe that
in general positive $p$-deficiency is a much stronger property than
large subgroup growth. In fact, if $\Gamma$ is a group, and $G$ is a
finite $p$-group with many generators, then the free product
$\Gamma\ast G$ should have smaller $p$-deficiency then $\Gamma$, while
most other asymptotic parameters should increase. More specifically,
we pose the following.

\begin{Problem}
Let $G$ be an elementary abelian group of order 64, $F_2$ a free group
with two generators. Decide whether $G\ast F_2$ has positive $p$-deficiency.
\end{Problem}

We believe that the Golod-Shafarevich inequality should
essentially hold true with deficiency replaced by $p$-deficiency. More
precisely, we pose the following.

\begin{Problem}
Let $G$ be a finite group. Prove that there exists a function $f(n)$, tending
to infinity with $n$, such that $\de(G)\leq - f(d(G))$.
\end{Problem}

It appears that the best way to attack this problem is by solving
Problem 1 first.

The concept of a $p$-Euler characteristic might also be generalized.

\begin{Problem}
Find a natural function, which is multiplicative on subgroups, and coincides
with the Euler-characteristic for virtually free groups, the
hyperbolic volume for Fuchsian groups, and $\chi_p$ for pro-$p$-groups.
\end{Problem}

Finally, one should be able to generalize the concept to discrete
groups without singling out a special prime $p$. However, the most obvious way
of doing so fails, for if $p, q$ are different prime numbers, then the
group $\langle x, y|x^{p^n}=y^{p^n}=x^{q^n}=y^{q^n}=1\rangle$ has four
relators, each of which only slightly affects the size of the group,
as, for example, measured by its subgroup growth. However, the
coprimality of the exponents forces the whole group to be trivial. We
believe that divisibility issues are the only cause for this kind of
collapse, more precisely, we pose the following.

\begin{Problem}
Let $w_1, \ldots$ an enumeration of all elements in $F_2$, and let
$(a_n)$ be a sequence of positive integers. Then the group $\Gamma$ as
$\Gamma=\langle x, y|w_1(x, y)^{a_1!}, w_2(x, y)^{a_2!}, \ldots\rangle$ is
obviously torsion. Prove that for every $\epsilon$ we can choose the
sequence $(a_n)$ in such a way that $s_n(\Gamma)\gg (n!)^{1-\epsilon}$.
\end{Problem}

\end{document}